\documentclass[letterpaper,12pt]{article}
\usepackage{fullpage} %
\usepackage{parskip} 
\usepackage{subfig}
\usepackage{tikz} 

\usepackage{mathtools}
\usepackage{bm}
\usepackage{physics}
\usepackage{amsfonts}
\usepackage{fancyhdr}
\usepackage{times}
\usepackage{changepage}
\usepackage{amssymb}
\usepackage{amsthm}
\usepackage[english]{babel}
\usepackage{graphicx}
\usepackage{xcolor,colortbl}
\usepackage{amsmath}
\usepackage{lipsum}
\usepackage{wrapfig}
\usepackage{float}
\usepackage{color} %red, green, blue, yellow, cyan, magenta, black, white
\usepackage{xcolor}
\usepackage{bbold}

\title{A Characterization of Distance Matrices of Positive Weighted Kneser Graphs and Generalized Petersen Graphs
}

\author{Joshua Steier, Luis Monterroso }

\newtheorem{theorem}{Theorem}

\begin{document}

\maketitle

\begin{abstract}
    
   Rubei et. al., established results for the distance matrix of positive weighted Petersen graphs. Focusing on the properties of the distance matrix, we generalized positive weighted Petersen graphs results to Kneser graphs. We analyzed theorems established by Rubei et al. and used girth of the generalized Petersen graphs and Kneser graphs to conclude generalizations. Further, we examined the properties of positive weighted generalized Petersen graphs. We generalized the properties of distance matrices of positive weighted Petersen graphs to positive weighted generalized Petersen graphs.

%needs some work


\end{abstract}
%to be completed 1/27/19
\section{Introduction:}
%rough sketch, make flow better
Graphs have a plethora of practical applications, for instance in network theory, and ecological dynamics. For instance, in [1], competitive hierarchy is described as a graph theoretic application to ecology.

A graph $G$ is an ordered pair $G=(V,E)$ composed of a set $V$ of vertices or nodes, together with a set $E$ of edges, which are 2-element subsets of $V$. A simple graph is a graph that does not contain multiple edges or loops. A weighted graph is a graph $G$ in which there is  a weight function $w(G):E(G)\rightarrow R^{+}$. In other words, a weighted graph contains a set of weights on each edge in the graph $G$.
In this manner, we may define a positive weight graph, as one in which the weights assigned to each edge are positive. The $k$-weight of $k$-subset of vertices $(v_{1},....v_{k})$ is the minimum among the weights of $G$ connecting $v_{1},...v_k$, which is the subgraphs of $G$ whose vertex set includes $v_{1}...v_{k}, $ denoted $D_{v1,...vk}(G)$. 

The Petersen graph is a special type of Kneser graph, which has been studied extensively and is used often as a counterexample.  In fact, Donald Knuth stated that the Petersen graph is "a remarkable configuration that serves as a counterexample to many optimistic predictions about what might be true for graphs in general."[2]  A Petersen graph contains 10 vertices and 15 edges, where it has one vertex for each 2-element subset of a 5-element set, and two vertices are connected by an edge if and only if the corresponding 2-element subsets are disjoint from each other.  The Petersen graph is a specific case of the Kneser graph, $KG_{5,2}$ since the Kneser graph itself is the graph whose vertices correspond to the $k$-element subsets of a set of $n$ elements, and where two vertices are adjacent if and only if the two corresponding sets are disjoint.

An odd graph is a graph that takes the form of $KG_{2n-1, n-1}$, and typically has high odd girth. The Petersen graph $KG_{5,2}$is isomorphic to the odd graph $O_3$.
Odd graphs have been used extensively in chemical graph theory and used as a generalization of the Petersen graph.

There also is a generalized Petersen graph, to which we feel is better suited for characterization from the Petersen graph $KG_{5,2}$. A generalized Petersen graph is denoted $G_{n,k}$ by Watkin's notation[3],  are a family of cubic graphs formed by connecting the vertices of a regular polygon to the corresponding vertices of a star polygon.
In Watkins' notation, $G(n,k)$ is a graph with vertex set
${u_0, u_1, ..., u_{n-1}, v_0, v_1, ..., v_{n-1}}$

and edge set
$u_i, u_{i+1}, {u_i v_i, v_i v_{i+ k}}: i = 0, ..., n-1$[3]

%notion of usefulness and indecomposable 

From previous work by Rubei et al[4], the definition of useful is In a positive-weighted graph $G$ an edge e is called useful if there exists
at least one couple of vertices $i$ and $j$ such that all the paths realizing $D_{i,j} (G)$ contain
the edge e. Otherwise, the edge is called useless.

%indecomposeable

Also from Rubei et al, the notion of indecomposable entry in the distance matrix is Let $D$ be a distance $m \cross m$ matrix for some $m \in N^{+}
. $We say that an
entry $D_{i,j}$ for some $i, j$ is indecomposable if and only if $D_{i,j} < D_{i,k} + D_{k,j}, 
\forall k ∈ {1, . . . , m} \ {i, j}.$

%good just finesse

\begin{section}{Findings and Methodology:}

The authors Rubei et al., sought to generalize properties of distance matrices of positive weighted Petersen graphs to Kneser graphs. These theorems are generalized results for Kneser graphs.

We define \#$X_1(x)$ as the set of the elements $y\in X$ such that $D_{x,y}$ is indecomposeable.

%proposed counterexample

\subsection{Distance matrices of Kneser graphs }

\begin{theorem}

Let $D$ be an $n \cross n$ distance matrix. Let us denote the set ${1,2,....,n}$ by $X$. 

The matrix $D$ is the distance matrix of a positively weighted Kneser graph in which each edge is useful if and only if the following hold:

a). for any $x\in X$, we have $\#X_{1}(x)= \binom{n-k}{k}$

%still determining if valid or invalid

%prove the other way

Since the girth for Kneser graphs can be 3 or 4, this means it is possible to have a cycle of length 3 or 4 in Kneser graphs. If we restrict Kneser graphs to those that are also odd graphs, and implement the following condition: $n\geq 3$, then the following will hold:

b). For any $k \in {3, 4},$ given $i_i, ..., i_k ∈ X$
such that $D_{i_1,i_2}, D_{i_2, i_3}, ..., D_{i_{t−1},i_t}$, are all indecomposable, we have that $D_{i_1,i_t}$ is decomposable. This can be generalized to subsets of $n-1$ by definition of the odd graph.

c). there exist distinct $v_{1},...,v_{n}\in X$ such that $D_{v_1, v_2, ...,}D_{v_4, v_n}, D_{v_n, v_1}$ are indecomposable and denoting $\hat{v_{j}}$ as the unique element in $X_1(v_j) - {v_1, v_2, v_3,...,v_n}$ we have that $\hat{v_1},...,\hat{v_n}$ are distinct.

d). 

$\forall x, y \in X$ with $y\not\in X_1(x)$ we have: 

$D_{x,y}= \min_{\substack{i_1,...,i_k\in X, k\geq 3, \\ i_1 = x, i_k=y, i_{j+1} \in X_1(i_j)\forall j}}$ $(D_{i_1, i_2}... + D_{i_{k-1}, i_k})$

Proof: Using proof techniques from [4],

a). Follows since Kneser graphs are, $\binom{n-k}{k}$-regular, so every vertex has degree $\binom{n-k}{k}$. b). Follows from the fact that in Kneser graphs , there are no cycles of length 3 or 4. Statement (c), can be seen from Figure 2, which can be extended to an arbitrary number of vertices $(v_1...v_n)$ in the odd graph. d). follows by definition of 2-weights, and the assumption that all edges are useful.

($\Leftarrow$) Let $v_1, ...,v_n, \hat{v_1}, ..., \hat{v_n}$ as in (c). Let $G$ be the graph in Figure 2 and for any adjacent vertices $i,j$ define $w(e(i,j))=D_{i,j}$ Let \textbf{G}= $(G, w).$ We want to show that $D_{x,y}(\textbf{G})= D_{x,y}$ for any $x,y \in X.$

Case 1: $D_{x,y}$ is indecomposeable.

First observe that, by (b), we have that $D_{vi, vj}$ is indecomposeable if and only if $(i,j)$ is one of the following: $(1,2), (2, 3), (3, 4), (4, 5), (5, 6), (7, 1)...$ This is the case for Kneser(7,3), which is also $O_4$ and can be readily extended for an arbitrary Kneser graph . We observe that $D_{vi,vj}$ is indecomposeable if and only if in the graph $G$ we have constructed there is an edge with endpoints $v_{i}$ and $v_{j}$.

Further, $D_{\hat{vi}, \hat{vj}}$ is indecomposeable if and only if in $G$ there is an edge with endpoints $\hat{v_i}$ and $\hat{v_j}$; otherwise there would be a contradiction with (b). 

Additionally, $D_{\hat{vi}, vj}$ is indecomposeable if and only if in $G$ there is an edge with endpoints $\hat{v_i}$ and ${v_j}$, otherwise this would contradict with b. Therefore, we can conclude that $D_{x,y}$ is indecomposeable if and only if in the graph $G$ we have constructed there is an edge with endpoints $x$ and $y$. In this case, $D_{x,y}(\textbf{G}= D_{x,y}$ by triangle inequalities.

Case 2:

$D_{x,y}$ is decomposable. By definition of 2-weights we have that $D_{x,y}(\textbf{G})$ is equal to number in (2). So, by the definition of $w$, it is equal to the number in (3). By the fact, we have proved before, that $D_{i, j}$ is indecomposeable if and only if in $G$ there is an edge with endpoints $i$ and $j$, by the decomposeability of $D_{x,y}$ and , finally by (d), we get that $D_{x,y}(\textbf{G})= D_{x,y}$

%must do other way

\end{theorem}

\subsection{Distance matrices of generalized Petersen graphs }

\begin{theorem}

Let $D$ be an $n \cross n$ distance matrix. Let us denote the set ${1,2,....,n}$ by $X$.

The matrix $D$ is the distance matrix of a positively weighted generalized Petersen graph in which each edge is useful if and only if the following hold:

a). for any $x\in X$, we have $\#X_{1}(x)= 3$

b).
So long as the following hold: $n\neq 3k$, $n \neq 4k$, and $k \neq 1$,([5]) for any $t \in {3, 4},$ given $i_i, ..., i_t ∈ X$
such that $D_{i_1,i_2}, D_{i_2, i_3}, ..., D_{i_{t−1},i_t}$, are all indecomposable, we have that $D_{i_1,i_t}$ is decomposable.

c). there exist distinct $v_{1},...,v_{n}\in X$ such that $D_{v_1, v_2, ...,}D_{v_4, v_n}, D_{v_n, v_1}$ are indecomposable and denoting $\hat{v_{j}}$ as the unique element in $X_1(v_j) - {v_1, v_2, v_3,...,v_n}$ we have that $\hat{v_1},...,\hat{v_n}$ are distinct.

d).
$\forall x, y \in X$ with $y\not\in X_1(x)$ we have: 

$D_{x,y}= \min_{\substack{i_1,...,ik\in X, k\geq 3, \\ i_1 = x, ik=y, i_{j+1} \in X_1(i_j)\forall j}}$ $(D_{i1, i2}... + D_{ik-1, ik})$

Proof:
Using proof techniques from [4]:

($\Rightarrow$) Assume that all the edges are useful, and an edge $e(i,j)$ is useful if and only if $D_{i,j}(G)$ is indecomposeable; so the edges correspond to indecomposable 2-weights.
a). Follows since generalized Petersen graphs are cubic, meaning they are 3-regular, so every vertex has degree 3. b). follows from the fact that in generalized Petersen graphs for the conditions described, there are no cycles of length 3 or 4. Statement (c), can be seen from Figure 1, which can be extended to an arbitrary number of vertices $(v_1...v_n)$ in the generalized Petersen graph. d). follows by definition of 2-weights, and the assumption that all edges are useful.

%prove the other way

($\Leftarrow$) Let $v_1, ...,v_n, \hat{v_1}, ..., \hat{v_n}$ as in (c). Let $G$ be the graph in Figure 1b and for any adjacent vertices $i,j$ define $w(e(i,j))=D_{i,j}$ Let \textbf{G}= $(G, w).$ We want to show that $D_{x,y}(\textbf{G})= D_{x,y}$ for any $x,y \in X.$

Case 1: $D_{x,y}$ is indecomposeable.

First observe that, by (b), we have that $D_{vi, vj}$ is indecomposeable if and only if $(i,j)$ is one of the following: $(1,2), (2, 3), (3, 4), (4, 5), (5,1).$ This is the case for generalized Petersen(5,2) and can be readily extended for an arbitrary generalized Petersen graph. We observe that $D_{vi,vj}$ is indecomposeable if and only if in the graph $G$ we have constructed there is an edge with endpoints $v_{i}$ and $v_{j}$.

Further, $D_{\hat{vi}, \hat{vj}}$ is indecomposeable if and only if in $G$ there is an edge with endpoints $\hat{v_i}$ and $\hat{v_j}$; otherwise there would be a contradiction with (b). 

Additionally, $D_{\hat{vi}, vj}$ is indecomposeable if and only if in $G$ there is an edge with endpoints $\hat{v_i}$ and ${v_j}$, otherwise this would contradict with b. Therefore, we can conclude that $D_{x,y}$ is indecomposeable if and only if in the graph $G$ we have constructed there is an edge with endpoints $x$ and $y$. In this case, $D_{x,y}(\textbf{G}= D_{x,y}$ by triangle inequalities.

Case 2:

$D_{x,y}$ is decomposable. By definition of 2-weights we have that $D_{x,y}(\textbf{G})$ is equal to number in (2). So, by the definition of $w$, it is equal to the number in (3). By the fact, we have proved before, that $D_{i, j}$ is indecomposeable if and only if in $G$ there is an edge with endpoints $i$ and $j$, by the decomposeability of $D_{x,y}$ and , finally by (d), we get that $D_{x,y}(\textbf{G})= D_{x,y}$

%must do other way

\end{theorem}

\end{section}

\begin{section}{Figures}

\begin{figure}[H]%
    \centering
    \subfloat[Generalized Petersen$(5,2)$]{{\includegraphics[width=5cm]{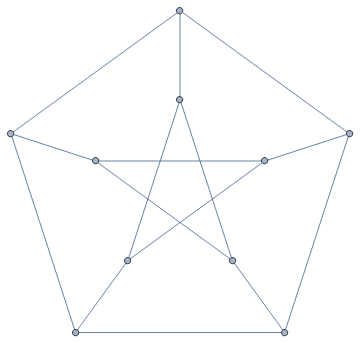} }}%
    \qquad
    \subfloat[Generalized Petersen $(10,2)$]{{\includegraphics[width=5cm]{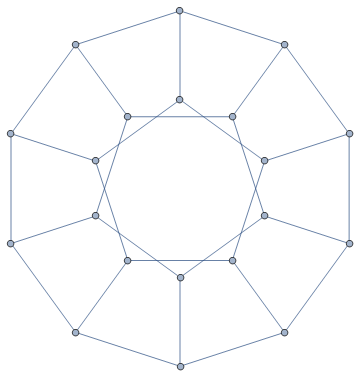} }}%
    \caption{Generalized Petersen $(5,2)$ and Generalized Petersen $(10,2)$}%
    \label{fig:example}%
\end{figure}

\begin{figure}[H]
    \centering
    \includegraphics[width=5cm]{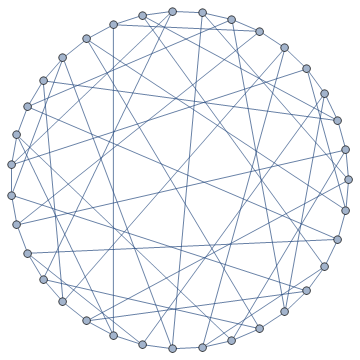}
    \caption{Kneser graph $(7,3)$}
    \label{fig:my_label}
\end{figure}
\end{section}

\begin{section}{Future Work}

In the future, there may be even more generalizations of the results obtained here. For instance, it may be possible that I graphs contain similar results. I graphs are generalizations of generalized Petersen graphs.

\end{section}

\begin{section}{References}

[1]	Dale, Mark, Applying Graph Theory in Ecological Research, Cambridge 2017

[2] Knuth, Donald E., The Art of Computer Programming; volume 4, pre-fasicicle 0A. A draft of section 7: Introduction to combinatorial searching

[3]  Watkins, Mark E.A Theorem on Tait colorings with an application to the generalized Petersen graphs. J. Combinational Theory 6 (1969), 152–164

[4] Rubei et al., A characterization of distance matrices of weighted
cubic graphs and Peterson graphs, arXiv, Jan. 2019
%how to cite? 

[5] Ferrero, Daniela; Hanusch, Sarah (2014), "Component connectivity of generalized Petersen graphs" (PDF), International Journal of Computer Mathematics, 91 (9): 1940–1963,

\end{section}

\end{document}